\documentclass[12pt,reqno]{amsart}
\usepackage{graphicx, amssymb, color,pdfsync}
\usepackage[all,cmtip]{xy}
\numberwithin{equation}{section}
\usepackage{mathrsfs}
\usepackage{hyperref}

\hypersetup{
    colorlinks=true,       % false: boxed links; true: colored links
    linkcolor=blue,          % color of internal links
    citecolor=blue,        % color of links to bibliography
    filecolor=blue,      % color of file links
    urlcolor=blue           % color of external links
}
\usepackage{tikz}
\usetikzlibrary{matrix,arrows}

\usepackage[switch]{lineno}
\usepackage{yfonts}
\usepackage{empheq} 
\advance\hoffset by -0.9 truecm

\begin{document}
\newcommand{\s}{\vspace{0.2cm}}
\newcommand{\Q}{\mathbb{Q}}

\newtheorem{theorem}{Theorem}
\newtheorem{prop}{Proposition}
\newtheorem{coro}{Corollary}
\newtheorem{lemm}{Lemma}
\newtheorem{claim}{Claim}
\newtheorem{example}{Example}
\theoremstyle{remark}
\newtheorem*{rema}{\it Remarks}
\newtheorem*{rema1}{\it Remark}
\newtheorem*{defi}{\it Definition}
\newtheorem*{theorem*}{\bf Theorem}
\newtheorem*{coro*}{Corollary}

\title[On the Jacobian of the Accola-Maclachlan curve of genus four]{On the Jacobian variety of the Accola-Maclachlan curve of genus four}
\date{}

\author{Robert Auffarth, Sebasti\'an Reyes-Carocca and Anita M. Rojas}
\address{Departamento de Matem\'aticas, Facultad de Ciencias, Universidad de Chile, Las Palmeras 3425, Santiago, Chile.}
\email{rauffar@uchile.cl, sebastianreyes.c@uchile.cl, anirojas@uchile.cl}

\thanks{This work was partially supported by Fondecyt Regular Grants 1220099, 1230708, 1230034 and 1220997, and MATH-AmSud Project 22-MATH-03}
%\keywords{Compact Riemann surfaces, group actions, automorphisms}
%\subjclass[2010]{30F10, 14H37, 30F35, 14H30}

\begin{abstract} In this short note, we study the Jacobian variety of the Accola-Maclachlan curve of genus four and obtain explicitly its Poincar\'e isogeny decomposition. More precisely, we show that its Jacobian variety is isomorphic to the product of two abelian surfaces that are simple, and provide explicitly a Riemann matrix for each one of the involved abelian surfaces.\end{abstract}
\maketitle

\section{Introduction and statement of the result}\label{sec:1}
\thispagestyle{empty}
A celebrated theorem due to Hurwitz, which goes back to the late nineteenth century \cite{H}, states that the automorphism group of a compact Riemann surface $X$ of genus $g \geqslant 2$ is finite, and that this order satisfies $$|\mbox{Aut}(X)| \leqslant 84(g-1).$$The smallest genus for which this upper bound is attained is three: the Klein's quartic $$x^3y+y^3z+z^3x=0$$represents a compact Riemann surface of genus three with 168 automorphisms. This famous curve was found by Klein himself a few years earlier than Hurwitz' Theorem was proved. It was not until seven decades later that MacBeath proved in \cite{M61} that there are infinitely many values of $g$ for which the Hurwitz' bound is attained.

\s

Let $N_g$ denote the maximum among  the orders of the automorphism groups of compact Riemann surfaces of a fixed genus $g \geqslant 2$. Since the Hurwitz Theorem provides a uniform upper bound for $N_g,$  the problem of determining a lower bound for $N_g$ arise naturally.  This problem was successfully solved independently by Accola in \cite{A68} and Maclachlan in \cite{M69}. More precisely, they showed that for each $g \geqslant 2$, the compact Riemann surface $X_g$ of genus $g$ represented by the algebraic curve  $$y^2=x^{2g+2}-1$$is endowed with a group of automorphisms of order $8(g+1);$ hence, they proved that $$N_g  \geqslant 8(g+1).$$Nowadays,  the curve $X_g$ is known as the {\it Accola-Maclachlan curve} of genus $g.$ 

\vspace{0.2 cm}

The role of $X_g$ became even more important because of the fact that for infinitely many values of $g$ the inequality above turns into an equality, showing that $8(g+1)$ is the best lower uniform bound for $N_g.$ In the early nineties, the problem of the uniqueness of $X_g$ was addressed by Kulkarni in \cite{K91}. Indeed, he proved that if $g \not\equiv 1 \mbox{ mod }4$ and $g$ is sufficiently large, then $X_g$ is the unique compact Riemann surface of genus $g$ endowed a group of automorphisms of order $8(g+1).$  

\vspace{0.2 cm}

We denote by $J_X$ the Jacobian variety of a compact Riemann surface $X$. The rule $$X \mapsto J_X \,\,\mbox{ defines a map }\,\, t : \mathscr{M}_g \to \mathscr{A}_g$$from the moduli space of compact Riemann surfaces of genus $g$ into the moduli space of  principally polarised abelian varieties of dimension $g.$ Torelli's theorem states that $t$ is injective, namely,  each isomorphism between  Jacobian varieties (as principally polarised abelian varieties) comes from an isomorphism between the corresponding Riemann surfaces. In other words, a Riemann surface is determined by its Jacobian variety: $$J_X \cong J_{X'} \,\, \iff \,\, X \cong X'.$$

Let $\mathscr{H}_g$ denote the Siegel upper space of $g \times g$ matrices; that is, the space of $g\times g$ complex symmetric matrices with positive definite imaginary part. The fact that $$\mathscr{A}_g \mbox{ is isomorphic to } \mathscr{H}_g/\mbox{Sp}(2g, \mathbb{Z}),$$where $\mbox{Sp}(2g, \mathbb{Z})$ denotes the symplectic group, shows that each principally polarised abelian variety $A$  can be described in terms of (the equivalence class of) a matrix $$Z_A \in \mathscr{H}_g.$$ Such a matrix is called a {\it Riemann matrix} for $A$.  However,  the problem of determining explicitly a Riemann matrix for a given abelian variety is very difficult, and has been successfully solved only for Jacobians of a few famous  Riemann surfaces. 

\vspace{0.2 cm}

By using the theory of Fuchsian groups, a general --but implicit-- form of Riemann matrices for the Jacobian variety of the Accola-Maclachlan curves was determined by Bujalance, Costa, Gamboa and Riera in \cite{BCGR00}.  They also used their results in order to provide an explicit Riemann matrix for $J_{X_2}$. See also \cite[\S.3]{M98}. Recently, an explicit Riemann matrix for $J_{X_4}$ was provided in \cite{RCR22} by considering a symplectic representation of the automorphism group of $X_4$ and supported by computational routines  developed in \cite{BRR13}.

 \vspace{0.2 cm}
  
We recall that an abelian variety is termed {\it simple} if it does not contain proper non-trivial abelian subvarieties.  Poincar\'e's irreducibility theorem states that if $A$ is an abelian variety then there are positive integers $n_j$ and simple abelian subvarieties $A_j$ of $A$ such that  $$A \sim A_1^{n_1} \times \cdots \times A_r^{n_r}$$where $\sim$ stands for isogeny. The decomposition above is called the {\it Poincar\'e isogeny decomposition} of $A$. It is worthwhile recalling here that both the problem of deciding whether or not a given abelian variety is simple, and the determination of its Poincar\'e isogeny decomposition are difficult tasks.

 \vspace{0.2 cm}

As we shall discuss later, the Jacobian variety $J_{X_g}$ is known to be non-simple and to have complex multiplication. In addition, for the case of genus $g=2$ one has that $$J_{X_2} \sim E_{-\omega} \times E_{\omega^2}$$where $\omega$ is a third root of unity. Observe that the Riemann matrices $(-\omega)$ and $(\omega^2)$ for the involved elliptic curves are defined over $\bar{\mathbb{Q}}$ as being elliptic curves with complex multiplication. This fact, together with a result due to Beauville \cite{B14}, implies that $$J_{X_2} \mbox{ is isomorphic to a product } E \times E'$$where $E\sim E'$ are  elliptic curves. Such curves can be chosen as $E_{i\sqrt{3}/3}$ and $E_{i\sqrt{3}}$.

\vspace{0.2 cm}

This short note is devoted to generalise the aforementioned facts, and to prove the following theorem.

\begin{theorem*} Let $J_{X_4}$ denote the Jacobian variety of the Accola-Maclachlan curve of genus four. There is an isomorphism of complex tori $$J_{X_4} \cong S_1 \times S_2$$where $S_1$ and $S_2$ are isogenous abelian subvarieties of $J_{X_4}$ of dimension two that are simple. If $\zeta$ denotes a fifth root of unity, then Riemann matrices for $S_1$ and $S_2$ are{\tiny$$Z_1=\left(\begin{array}{rr}
-8 \zeta^{3} + 4 \zeta^{2} - 4 \zeta - 2 & -4 \zeta^{3} - 4 \zeta^{2} - 6 \\
-4 \zeta^{3} - 4 \zeta^{2} - 6 & 4 \zeta^{3} + 4 \zeta + 2
\end{array}\right) \,\,\,\, Z_2=\left(\begin{array}{rr}
-\frac{4}{5} \zeta^{3} + \frac{8}{5} \zeta^{2} + \frac{4}{5} \zeta + \frac{2}{5} & -\frac{4}{5} \zeta^{2} - \frac{4}{5} \zeta - \frac{2}{5} \\
-\frac{4}{5} \zeta^{2} - \frac{4}{5} \zeta - \frac{2}{5} & \frac{4}{5} \zeta^{3} + \frac{4}{5} \zeta^{2} + \frac{8}{5} \zeta + \frac{4}{5}
\end{array}\right)$$ }

%the matrix
%
%
%
% $${\tiny\left(\begin{array}{rr}
%-4 \zeta^{3} + 2 \zeta^{2} - 2 \zeta - 1 & -2 \zeta^{3} - 2 \zeta^{2} - 3 \\
%-2 \zeta^{3} - 2 \zeta^{2} - 3 & 2 \zeta^{3} + 2 \zeta + 1
%\end{array}\right)}$$is a Riemann matrix for $S_1$ and $S_2$. 
\end{theorem*}

Observe that the theorem above provides explicitly the Poincar\'e isogeny decomposition of $J_{X_4}.$ At the end of this paper we shall briefly discuss the situation for $J_{X_6}$  as well as some open questions for the general case.

 \section{Preliminaries}\label{sec:2}

\subsection*{The Accola-Maclachlan curve of genus $g$}

Let $g \geqslant 2$ be an integer and let $T \subset \mathbb{H}$ be a hyperbolic triangle of angles $\tfrac{\pi}{2}$, $\tfrac{\pi}{4}, \tfrac{\pi}{2g+2}.$ If $R_1, R_2$ and $R_3$ are the reflexions in each one of the sides  of $T$ then $$\gamma_1=R_3\circ R_1, \,\, \gamma_2=R_1 \circ R_2 \,\, \mbox{ and }\,\, \gamma_3=R_2 \circ R_3$$generate a discrete group $\Gamma$ of automorphisms of  $\mathbb{H}.$ It follows that the quotient $\mathbb{H}/\Gamma$ has a structure of (compact) Riemann orbifold. The group $\Gamma$, usually called a {\it triangle Fuchsian group of signature $(0; 2,4,2g+2)$}, has a canonical presentation $$\Gamma=\langle \gamma_1, \gamma_2, \gamma_3 : \gamma_1^2=\gamma_2^4=\gamma_3^{2g+2}=\gamma_1\gamma_2\gamma_3=1\rangle.$$In addition, the elements of $\Gamma$ of finite order are the elements conjugate to powers of $\gamma_1, \gamma_2$ and $\gamma_3.$ We now consider the finite group $G$ of order $8(g+1)$ with presentation$$G=\langle a,b,c: a^{2g+2}=b^2=c^2= [a,b]=[b,c]=cacab=1 \rangle$$ (a semidirect product of $\mathbb{Z}_{2g+2} \times \mathbb{Z}_2$ and $\mathbb{Z}_2).$ The correspondence $$\gamma_1 \mapsto c, \,\, \gamma_2 \mapsto ca, \,\, \gamma_3 \mapsto a^{-1}$$ defines a epimorphism $\Gamma \to G$ whose kernel $K$ does not have nontrivial elements of finite order. Such an epimorphism is usually called a {\it surface-kernel epimorphism}. It then follows from the uniformisation theorem that the  compact Riemann surface  $$X:=\mathbb{H}/K$$ is endowed with a group of automorphisms $\Gamma/K$ isomorphic to  $G.$ Moreover, the associated branched regular covering map induced by the inclusion $K \leqslant \Gamma$$$X=\mathbb{H}/K \to \mathbb{H}/\Gamma \cong G \cong \mathbb{P}^1$$ramifies over exactly three values,  marked with $2, 4$ and $2g+2.$ It follows from the Riemann-Hurwitz formula that the genus of $X$ is $g.$ Note that $X$ has a central involution, represented by $b \in G,$ fixing exactly $2g+2$ points. Hence, $X$ is hyperelliptic or, equivalently,   $X$ is isomorphic to the normalisation of an algebraic  curve of the form $$y^2=(x-t_1)\cdots (x-t_{2g+2})$$where $t_1, \ldots, t_{2g+2}$ are pairwise distinct complex numbers. In addition, since $G/\langle y \rangle$ is isomorphic to a dihedral group, after a suitable M\"{o}bius transformation, the complex numbers $t_1, \ldots, t_{2g+2}$ can be taken as the $(2g+2)$-th roots of unity. Consequently, $X$ is isomorphic to the Accola-Maclachlan curve $X_g$. It follows from classical results due to Singerman \cite{S72} that the (full) automorphism group of $X_g$ is isomorphic to $G.$

\subsection*{Riemann matrices and abelian varieties}
Let $T=V/\Lambda$ be a complex torus of dimension $g$. If we choose bases $$B=\{v_1, \ldots, v_g\} \,\, \mbox{ and } \,\,  L=\{\lambda_1, \ldots, \lambda_{2g}\}$$ of $V$ and $\Lambda$ as a complex vector space and  as a $\mathbb{Z}$-module respectively, then the  {\it period matrix} for $T$ with respect to these bases is the $g \times 2g$ matrix defined as $$\Pi = (\pi_{i,j}) \mbox{ where } \, \lambda_{j}= \pi_{1,1}v_1+\cdots + \pi_{1,g}v_g \mbox{ for each }j.$$Moreover, by choosing appropriate bases, the period matrix has the form $$\Pi=(\mbox{I}_g \, Z) \mbox{ where } \det(\mbox{Im}(Z)) \neq 0,$$with $\mbox{I}_g$ denoting the identity $g \times g$ matrix. 

\vspace{0.2 cm}

Let $f: T \to T'$ be a homomorphism of complex tori of dimension $g$ and $g'$, and let $\Pi$ and $\Pi'$ be their period matrices with respect to given bases. One can associate to $f$ two linear maps: the {\it analytic representation} and the {\it  rational representation}  $$\rho_a(f) : V \to V' \,\, \mbox{ and } \,\, \rho_r(f) : \Lambda \to \Lambda',$$represented by a $g' \times g$ complex matrix and by a $2g' \times 2g$ integral matrix respectively. The fundamental relation they satisfy is given by
\begin{equation*} \label{eq:Hurwitz}
	\rho_a(f) \ \Pi = \Pi' \, \rho_r(f).
\end{equation*} 

Observe that if $f$ is an isomorphism then $g'=g$, the matrix $\rho_a(f)$ is  nonsingular, and the matrix $\rho_r(f)$ is unimodular. 
\vspace{0.2 cm}

A complex {\it abelian variety} is a complex torus which is also a complex projective algebraic variety. Each abelian variety $T=V/\Lambda$ admits a {\it polarisation},  that is, a non-degenerate real alternating form $\Theta$ on $V$ such that for all $v,w \in V$$$\Theta(iv, iw)=\Theta(v,w) \,\,  \mbox{ and } \,\, \Theta(\Lambda \times \Lambda) \subset \mathbb{Z}.$$

If each elementary divisor of $\Theta|_{\Lambda \times \Lambda}$ equals 1 then
$\Theta$ is called {\it principal} and  $T$ is called a
 {\it principally polarised abelian variety}; from now on, we write  {\it ppav} for short. In this case,  there exists a basis for $\Lambda$ such that the matrix for $\Theta_{\Lambda \times \Lambda}$  is given by \begin{equation*}\label{simpl}
J = \left( \begin{smallmatrix}
0 & \mbox{I}_g \\
-\mbox{I}_g & 0
\end{smallmatrix} \right) \,\, \mbox{ where } \,\, g= \dim(X);
\end{equation*}such a basis  is called {\it symplectic}. In addition, there exist a basis for $V$ with respect to which the period matrix for $T$ is $\Pi=(\mbox{I}_g \, Z)$ where $$Z  \in \mathscr{H}_g=\{ Z \in \mbox{M}(g, \mathbb{C}) : Z = Z^t , \, \mbox{Im}(Z) >0\},$$with $Z^t$ denoting the transpose matrix of $Z.$ By an {\it isomorphism of ppavs} we mean an isomorphism of the underlying complex tori  preserving the involved polarisations. Observe that if $(\mbox{I}_{g} \, Z_i)$ is the period matrix for $T_i$ then an isomorphism $T_1 \to T_2$ is given by invertible matrices \begin{equation}\label{ig}M \in \mbox{GL}(g, \mathbb{C}) \,\, \mbox{ and } \,\, R \in \mbox{GL}(2g, \mathbb{Z}) \,\, \mbox{ such that } \,\,  M(\mbox{I}_{g} \, Z_1)=(\mbox{I}_{g} \, Z_2)R.\end{equation}

As $R$ preserves the polarisation $J$, it must satisfy $R^tJ R=J $ and hence it belongs to $\mbox{Sp}(2g, \mathbb{Z})$. It follows from \eqref{ig} that the orbits of the action $\mbox{Sp}(2g, \mathbb{Z}) \times \mathscr{H}_g \to \mathscr{H}_g$ given by
\begin{equation*}       (R=  \left( \begin{smallmatrix}
A & B \\
C & D
\end{smallmatrix} \right) , Z ) \mapsto R \cdot Z := (A+ZC)^{-1}(B+ZD)
\end{equation*}are in bijective correspondence with the isomorphism classes of ppavs of dimension $g.$ Hence $$\mathscr{H}_g \to \mathcal{A}_g:=\mathscr{H}_g/ \mbox{Sp}(2g, \mathbb{Z})$$is the moduli space of isomorphism classes of ppavs of dimension $g.$ See \cite{BL} and \cite{O76} for more details.

\section{Proof of the theorem}

{\bf The decomposition of $J_{X_4}$.} We recall that the automorphism group of the  Accola-Maclachlan curve $X_4$ of genus four is isomorphic to 
$$G=\langle a,b,c: a^{10}=b^2=c^2= [a,b]=[b,c]=cacab=1 \rangle,$$and that the action of $G$ on $X_4$ is represented by the surface-kernel epimorphism 
$$\theta=(c,ca^{-1}, a)$$

We write $J_{X_4}=V/\Lambda$ where $V$ is the $g$-dimensional complex vector space of holomorphic $1$-forms on $X_4$ and $L$ is the first (integral) homology group of $X_4.$

\vspace{0.2 cm}

We follow the routines developed in \cite{BRR13} to ensure that there is a symplectic basis $\gamma=\{\alpha_1,\ldots, \gamma_4, \beta_1,\ldots,\beta_4\}$ for $L$  (and hence  $\{\alpha_1, \ldots, \alpha_4\}$ is a basis for $V$) in such a way that the rational representation  $$\rho_r:G\to \mbox{GL}(L\otimes \mathbb{Q})$$ induced by the action of $G$  is given by
$$ \rho_r(a)=\left(\begin{array}{rrrrrrrr}
0 & 1 & 1 & 1 & -1 & 0 & 0 & 0 \\
1 & 0 & 0 & 0 & 0 & -1 & 0 & 0 \\
0 & 0 & 0 & 0 & 0 & 1 & -1 & 0 \\
0 & 0 & 0 & 0 & 0 & 0 & 1 & -1 \\
1 & 0 & 0 & 0 & 0 & 0 & 0 & 0 \\
0 & 1 & 1 & 1 & 0 & 0 & 0 & 0 \\
0 & 0 & 1 & 1 & 0 & 0 & 0 & 0 \\
0 & 0 & 0 & 1 & 0 & 0 & 0 & 0
\end{array}\right)\,\,$$ $$ \rho_r(c)=\left(\begin{array}{rrrrrrrr}
0 & 0 & 0 & 0 & 0 & 0 & -1 & 1 \\
0 & -1 & 0 & 0 & 0 & 0 & 0 & 0 \\
0 & 0 & 0 & 0 & 1 & 0 & 0 & 0 \\
0 & 1 & 1 & 1 & -1 & 0 & 0 & 0 \\
0 & 0 & 1 & 0 & 0 & 0 & 0 & 0 \\
0 & 0 & 0 & 0 & 0 & -1 & 0 & 1 \\
-1 & 0 & 0 & 0 & 0 & 0 & 0 & 1 \\
0 & 0 & 0 & 0 & 0 & 0 & 0 & 1
\end{array}\right)$$(observe that $b=(ac)^2$). 

Although an explicit Riemann matrix for $J_{X_4}$ has been found in \cite{RCR22},  we prefer to compute it again here by finding the fixed points in $\mathscr{H}_4$ by the action of $$G_r=\langle \rho_r(a), \rho_s(c)\rangle \leqslant \mbox{Sp}(8, \mathbb{Z}).$$In other words, a period matrix for $J_{X_4}$ is the (unique) solution $Z$ of the equation $$R \cdot Z =R \mbox{ for all } R \in G_r.$$ Explicitly, if $\zeta$ a fifth root of unity then {\tiny$$Z=\left(\begin{array}{rrrr}
-2 \zeta^{3} + \zeta^{2} - \zeta - \frac{1}{2} & -\zeta^{3} - \zeta^{2} - 2 & \zeta^{3} + \zeta^{2} + \frac{3}{2} & \zeta^{3} + \zeta^{2} + \frac{3}{2} \\
-\zeta^{3} - \zeta^{2} - 2 & \frac{6}{5} \zeta^{3} + \frac{2}{5} \zeta^{2} + \frac{8}{5} \zeta + \frac{4}{5} & -\frac{3}{5} \zeta^{3} - \frac{1}{5} \zeta^{2} - \frac{4}{5} \zeta - \frac{2}{5} & -\frac{7}{5} \zeta^{3} + \frac{1}{5} \zeta^{2} - \frac{6}{5} \zeta - \frac{3}{5} \\
\zeta^{3} + \zeta^{2} + \frac{3}{2} & -\frac{3}{5} \zeta^{3} - \frac{1}{5} \zeta^{2} - \frac{4}{5} \zeta - \frac{2}{5} & \frac{4}{5} \zeta^{3} + \frac{3}{5} \zeta^{2} + \frac{7}{5} \zeta + \frac{7}{10} & \frac{6}{5} \zeta^{3} - \frac{3}{5} \zeta^{2} + \frac{3}{5} \zeta + \frac{3}{10} \\
\zeta^{3} + \zeta^{2} + \frac{3}{2} & -\frac{7}{5} \zeta^{3} + \frac{1}{5} \zeta^{2} - \frac{6}{5} \zeta - \frac{3}{5} & \frac{6}{5} \zeta^{3} - \frac{3}{5} \zeta^{2} + \frac{3}{5} \zeta + \frac{3}{10} & \frac{4}{5} \zeta^{3} + \frac{3}{5} \zeta^{2} + \frac{7}{5} \zeta + \frac{7}{10}
\end{array}\right).$$}

We point out that the period matrix $\Pi =(\mbox{I}_4 \, Z)$ we found for $J_{X_4}$ allows us to recover each $\beta_j$ in terms of a linear combination of $\alpha_1, \ldots, \alpha_4.$

\vspace{0.2 cm}
 
We consider the subgroup $H_1=\langle c\rangle$ of $G$. Following \cite[Theorem 3]{RCR22}, there is an isogeny $$J_{X_4} \to J_{Y_1}^2$$where $Y_1$ is the genus two quotient Riemann surface associated to the action of $H_1$ on $X_4$. In addition, if we consider the subgroup $H_2=H^a$ of $G$ and  accordingly we denote the corresponding quotient Riemann surface by $Y_2$, then $J_{Y_1}$ and  $J_{Y_2}$ are isomorphic (as ppavs). Hence,  there is an isogeny   $$\tilde{\varphi} : J_{Y_1} \times J_{Y_2} \to J_{X_4}. $$

We claim that $\tilde{\varphi}$ is an isomorphism of complex tori.

\vspace{0.2 cm}

In order to prove the claim, we start by recalling that the regular covering maps  $$\pi_j:X_4\to Y_j \,\, \mbox{ for }\,\, j=1,2 $$ ramify over exactly two values. It follows from \cite[Proposition 11.4.3]{BL} that the pull-back $\pi^*_j: J_{Y_j} \to J_{X_4}$ is injective and consequently it induces an isomorphic onto its image, namely, $$J_{Y_j} \cong S_j:=\pi_j^*(J_{Y_j}) \subset J_{X_4} \, \mbox{ for } \, j=1,2.$$The subvarieties $S_1$ and $S_2$ of $JX_4$ can be described in a very explicit way as follows. We denote by $$\Phi: \mathbb{Q}[G] \to \mbox{End}_{\mathbb{Q}}(J_{X_4})$$the $\mathbb{Q}$-algebra homomorphism induced by the action of $G$ on $X_4.$ Then  $$S_j:=\mbox{Im}(\Phi(p_{H_j})) \,\, \mbox{ where } \,\, p_{H_j}=\tfrac{1}{2}\Sigma_{h \in H_j}h \in \mathbb{Q}[G]$$(see \cite{cr} for more details). Further, we can describe $S_j=V_j/\Lambda_j$ in terms of coordinates as follows.  The rational representation of $p_{H_j}$ is $$\rho_r(p_{H_j})=\tfrac{1}{2}\Sigma_{h \in H_j}\rho_r(h) \in \mbox{Sp}(8, \mathbb{Q})$$ and the coordinates of $\Lambda_j$ with respect to  the symplectic basis $\gamma$ are given by$$B_1=\left(\begin{array}{rrrrrrrr}
0 & 0 & 0 & 1 & 0 & 0 & 0 & 0 \\
1 & 0 & 0 & 0 & 0 & 0 & -1 & 0 \\
1 & 0 & 0 & 0 & 0 & 1 & 1 & 2 \\
0 & 0 & 1 & 0 & 1 & 0 & 0 & 0
\end{array}\right)\,\,  \text{ and } \,\, B_{2}=\left(\begin{array}{rrrrrrrr}
0 & 0 & 0 & 1 & -1 & 0 & 0 & 0 \\
0 & 1 & 1 & 0 & 0 & 0 & 0 & 0 \\
1 & 0 & 0 & 0 & 0 & 0 & 0 & 1 \\
0 & 0 & 0 & 0 & 0 & 1 & 1 & 0
\end{array}\right)$$

Once a period matrix $\Pi$ for $J_{X_4}$ is given, the procedures developed in \cite{rr-preprint} (assisted by the algorithms in \cite{rrgithub}) allow us to identify the relations between the coordinates captured in the matrices $B_j$. We obtain that
\begin{enumerate}
\item the induced polarisation in the subvarieties $S_j$ is twice a principal one, and that

\s

\item Riemann matrices for $S_1$ and $S_2$ are{\tiny$$Z_1=\left(\begin{array}{rr}
-8 \zeta^{3} + 4 \zeta^{2} - 4 \zeta - 2 & -4 \zeta^{3} - 4 \zeta^{2} - 6 \\
-4 \zeta^{3} - 4 \zeta^{2} - 6 & 4 \zeta^{3} + 4 \zeta + 2
\end{array}\right) \,\,\,\, Z_2=\left(\begin{array}{rr}
-\frac{4}{5} \zeta^{3} + \frac{8}{5} \zeta^{2} + \frac{4}{5} \zeta + \frac{2}{5} & -\frac{4}{5} \zeta^{2} - \frac{4}{5} \zeta - \frac{2}{5} \\
-\frac{4}{5} \zeta^{2} - \frac{4}{5} \zeta - \frac{2}{5} & \frac{4}{5} \zeta^{3} + \frac{4}{5} \zeta^{2} + \frac{8}{5} \zeta + \frac{4}{5}
\end{array}\right)$$ }
\end{enumerate}The fact that each $S_j$ is an abelian subvariety of $J_{X_4}$ implies that the sum induces an isogeny $$\varphi: S_1\times S_2\to J_{X_4}.$$Observe that $\varphi \circ (\pi_1^*, \pi_2^*)=\tilde{\varphi}.$ The  rational representations of $\varphi$ is$$\rho_r(\varphi)=\left(\begin{array}{rrrrrrrr}
0 & 1 & 0 & 0 & 1 & 0 & 1 & 0 \\
0 & 0 & 0 & 1 & 0 & 0 & 0 & 0 \\
0 & 0 & 0 & 1 & 0 & 1 & 0 & 0 \\
1 & 0 & 1 & 0 & 0 & 0 & 0 & 0 \\
0 & 0 & -1 & 0 & 0 & 1 & 0 & 0 \\
0 & 0 & 0 & 0 & 1 & 0 & 0 & 1 \\
0 & -1 & 0 & 0 & 1 & 0 & 0 & 1 \\
0 & 0 & 0 & 0 & 2 & 0 & 1 & 0
\end{array}\right).$$Finally, since the determinant of the latter matrix equals one (and this number agrees with the order of $\mbox{ker}(\varphi)$), we conclude that $\varphi$ is an isomorphism.

\vspace{0.2 cm}

{\bf The simplicity of $S_1$ and $S_2$.} Since $S_1$ and $S_2$ are isomorphic ppavs, it is enough to show that $S_1$ is simple.  To accomplish this task, we employ the equivalence between the following statements, proved in \cite[Proposition 4.4]{ALR17}.

\begin{enumerate}
\item The abelian surface represented by the Riemann matrix $$\left(\begin{array}{cc}\tau_{1}&\tau_2\\\tau_2&\tau_3\end{array}\right)$$admits a 
sub elliptic curve. 
\item There exists a $6$-uple  
$(a_{12},a_{13},a_{14},a_{23},a_{24},a_{34})$ of rational numbers satisfying \begin{equation}-1=a_{13}+a_{24} \tag{i}\end{equation}\begin{equation}(\tau_1\tau_3-\tau_2^2)a_{12} - a_{14}\tau_1 + a_{13}\tau_2 - a_{24}\tau_2 + a_{23}\tau_3 + a_{34}=0\tag{ii}\end{equation} \begin{equation}a_{14}a_{23}-a_{13}a_{24}+a_{12}a_{34}=0. \tag{iii}\end{equation}
\end{enumerate}

We apply the criterion above with the Riemann matrix for $S_1$ $$\tilde{Z}=\tfrac{1}{2}Z_1={\tiny\left(\begin{array}{rr}
-4 \zeta^{3} + 2 \zeta^{2} - 2 \zeta - 1 & -2 \zeta^{3} - 2 \zeta^{2} - 3 \\
-2 \zeta^{3} - 2 \zeta^{2} - 3 & 2 \zeta^{3} + 2 \zeta + 1
\end{array}\right)}$$
which has determinant $\det(\tilde{Z})=2 \zeta^{3} + 2 \zeta^{2} + 2$. By replacing the entries of $\tilde{Z}$ in equation (ii), we obtain a linear combination over the rationals of the form $$A+B\zeta+C\zeta^2+D\zeta^3=0.$$Now, the fact that $\{1, \zeta,\zeta^2,\zeta^3\}$ is linearly independent over $\Q$ yields four equations, labeled as (3.1b) to (3.1e) below. The system turns into\begin{subequations}
  \begin{empheq}{align}
    a_{13} + a_{24} = -1 \\
    a_{14}+a_{23}=0\\ 
    a_{12}-a_{13}-a_{14}-a_{24}=0\\
     a_{12}-a_{13}+2a_{14}+a_{23}+a_{24}=0\\
     2a_{12}-3a_{13}+a_{14}+a_{23}+3a_{24}+a_{34}=0\\
     a_{14}a_{23}-a_{13}a_{24}+a_{12}a_{34}=0
  \end{empheq}
\end{subequations}Observe that equations (3.1a) to (3.1e) give rise to a linear system of equations in six variables, whose solutions are $$(a_{12},a_{13},a_{14},a_{23},a_{24},a_{34})=(\mu,\tfrac{\mu-1}{2}, 0,0,-\tfrac{1+\mu}{2}, \mu) \,\, \mbox{ where } \mu \in \mathbb{C}.$$We replace this information in equation (3.1f) to obtain that

$$\tfrac{\mu-1}{2} \cdot \tfrac{\mu+1}{2}+\mu^2=0 \,\, \mbox{ and therefore } \,\,  5\mu^2=1.$$Hence, the system does not have any rational solution, showing that $S_1$ is simple. Since $S_1$ is isogenous to $S_2$, we get that $S_2$ is also simple. 
Hence, we have the Poincaré decomposition of JX4.

\section{Final comments}

We end this article by mentioning some comments and remarks.

\vspace{0.2cm}

{\bf CM type and Riemann matrices.} It is known that $J_{X_g}$ admits complex multiplication, for each $g \geqslant 2$. For instance, this property follows from \cite[Theorem 4]{W00}, after noticing that the regular covering map $X_g \to X_g/\langle a \rangle \cong \mathbb{P}^1$ ramifies over three values. If we further assume that $g+1$ is prime then, as proved in  \cite{RCR22},  there is an isogeny between $J_{X_g}$ and $J_{Y}^2$, where $Y$ is the quotient Riemann surface associated to the action of $\langle c \rangle$ on $X_g$. Thus, $Y$  admits complex multiplication too. Moreover, by the results of \cite[\S4]{G05},  $Y$ is defined over $\bar{\mathbb{Q}}$. Now, since $J_Y$ admits complex multiplication and is defined over $\bar{\mathbb{Q}}$, a result due to Wolfart \cite{W00} ensures that each Riemann matrix for $J_Y$ is a $\bar{\mathbb{Q}}$-rational point of the Siegel upper half-space. For genus four we went further, and obtained that a Riemann matrix for $J_Y$ is $${\tiny\left(\begin{array}{rr}
-4 \zeta^{3} + 2 \zeta^{2} - 2 \zeta - 1 & -2 \zeta^{3} - 2 \zeta^{2} - 3 \\
-2 \zeta^{3} - 2 \zeta^{2} - 3 & 2 \zeta^{3} + 2 \zeta + 1
\end{array}\right)}$$which is, as expected, a $\bar{\mathbb{Q}}$-rational point of $\mathscr{H}_2.$

\vspace{0.2 cm}

{\bf The case of $JX_6$ and some open problems.} Assume that $g+1$ is a prime number. By arguing as done in the case of genus four, it can be seen that   \begin{equation}\label{hola}J_{Y_1} \times J_{Y_2} \cong A_1 \times A_2 \sim J_{X_g},\end{equation}where $A_1$ and $A_2$ are abelian subvarieties of $J_{X_g}$ of dimension $\tfrac{g}{2}.$ These subvarieties are, in general, not complementary. In fact, they are associated to pairs of subgroups $H_1$ and $H_2$ of $G$ conjugate to $\langle c \rangle$ and the corresponding idempotents $p_{H_1}$ and $p_{H_2}$ do not sum the identity. It would be interesting to understand the extent to which our methods can be generalised to the general case. For instance, in the case of genus six, by choosing appropriate subgroups $H_1$ and $H_2$ as before, it can be seen that the symplectic representation of the corresponding isogeny in \eqref{hola} has determinant one. Hence, the isogeny  is an isomorphisms.

In fact, more can be said about the case $g=6$. By calculating a period matrix for $A_1$ (which is isogenous to $A_2$) and calculating its N\'eron-Severi group, it can be seen that its Picard number is $9=3^2$. By the main theorem of \cite{K75}, we have that this implies that $A_1$ (and therefore $A_2$) is isogenous to the self-product of an elliptic curve (with complex multiplication). In particular, $JX_6$ itself is isogenous to the self-product of an elliptic curve (with complex multiplication),  as happens in the case $g=2$. This naturally leads us to the following amusing question that we would like to answer in the future:\\

\subsection*{Question:} How many primes $p$ exist such that $p\equiv 3 \mbox{ mod }4$ and $JX_{p-1}$ splits isogenously as the product of elliptic curves (with complex multiplication)?


\begin{thebibliography}{99}


\bibitem{A68}  
Accola, R. D. M., 
On the number of automorphisms of a closed Riemann surface. 
Trans. Amer. Math. Soc. {\bf 131} (1968), 398--408. 

\bibitem{ALR17} Auffarth, R.; Lange, H.; Rojas, A. M., 
A criterion for an abelian variety to be non-simple. 
J. Pure Appl. Algebra {\bf 221} (2017), no. 8, 1906--1925.

\bibitem{B14} 
Beauville, A., 
Some surfaces with maximal Picard number. 
J. \'{E}c. polytech. Math. {\bf 1} (2014), 101--116.
 
\bibitem{BRR13} 
Behn, A.; Rodríguez, R. E.; Rojas, A. M., 
Adapted hyperbolic polygons and symplectic representations for group actions on Riemann surfaces,  
J. Pure Appl. Algebra {\bf 217} (2013), no. 3, 409--426.
  
\bibitem{BL} 
Birkenhake, C.; Lange, H., 
Complex abelian varieties. Second edition. 
Grundlehren der mathematischen Wissenschaften {\bf 302}. Springer-Verlag, Berlin, 2004. xii+635 pp. ISBN: 3-540-20488-1.
 
\bibitem{BCGR00} 
Bujalance, E.; Costa, A. F.; Gamboa, J. M.; Riera, G. 
Period matrices of Accola-Maclachlan and Kulkarni surfaces. 
Ann. Acad. Sci. Fenn. Math. 25 (2000), no. 1, 161--177. 

\bibitem{cr}
Carocca, A.; Rodr\'iguez, R. E.; 
Jacobians with group actions and rational idempotents, 
J. Algebra {\bf 306} (2006), no. 2, 322--343.

\bibitem{G05} 
Gonz\'alez-Diez, G., 
Variations on Belyi's theorem. 
Q. J. Math. {\bf 57} (2006), no. 3, 339--354. 

\bibitem{H} Hurwitz, A. 
Ueber algebraische Gebilde mit eindeutigen Transformationen in sich. (German) 
Math. Ann. 41 (1892), no. 3, 403--442.


\bibitem{K75} Katsura, T., 
On the structure of singular abelian varieties, 
Proc. Japan Acad., 51(4):224-228, 1975.

\bibitem{K91} Kulkarni, R. S., 
A note on Wiman and Accola-Maclachlan surfaces. 
Ann. Acad. Sci. Fenn. Ser. A I Math. {\bf 16} (1991), no. 1, 83--94. 

\bibitem{M61} Macbeath, A. M., 
On a theorem of Hurwitz, Proc. 
Glasgow Math. Assoc. {\bf 5} (1961), 90--96.

\bibitem{M69} Maclachlan, C., 
A bound for the number of automorphisms of a compact Riemann surface. 
J. London Math. Soc. {\bf 44} (1969), 265--272. 

\bibitem{M98} Marchisio, M. R., 
Abelian surfaces and products of elliptic curves. 
Boll. Unione Mat. Ital. Sez. B Artic. Ric. Mat. (8) {\bf 1} (1998), no. 2, 407--427.

\bibitem{O76} Oort F., 
Singularities of coarse moduli schemes, 
S\'em. Dubriel {\bf 16} (1976).
 
\bibitem{RCR22} Reyes-Carocca, S.; Rojas, A. M., 
On large prime actions on Riemann surfaces. 
J. Group Theory 25 (2022), no. 5, 887--940. 

\bibitem{rr-preprint}
Rodr\'iguez, R. E; Rojas,  A. M.,
Period matrices for abelian varieties: 
Algorithms and applications. Preprint (2024).

\bibitem{rrgithub}
Rodr\'iguez, R. E; Rojas,  A. M.,
Period Matrices (Magma code). 
Available at \url{https://github.com/rojas-ani/Period-matrices}.

\bibitem{S72} 
Singerman, D., 
Finitely maximal Fuchsian groups. 
J. London Math. Soc. (2) 6 (1972), 29--38.
 
\bibitem{W00} 
Wolfart, J., Triangle groups and Jacobians of CM type. 
available at \url{https://www.math.uni-frankfurt.de/~wolfart/Artikel/jac.pdf}, Frankfurt a.M., 2000

\end{thebibliography}
\end{document}